\documentclass[12pt]{amsart}
\usepackage{amsmath}
\usepackage{amsfonts,amssymb,amsthm}

\newtheorem{theorem}{Theorem}[section]

\begin{document}


\title[A Remark]{A Remark on ``Counting primitive elements in free groups'' (by J.~Burillo and E.~Ventura)}

\author{Igor Rivin}


\address{Department of Mathematics, Temple University, Philadelphia}

\curraddr{Mathematics Department, Princeton University}

\email{rivin@math.temple.edu}

\thanks{The author was supported by a grant from the National Science Foundation.}

\date{today}

\keywords{primitive elements, free groups, punctured tori}

\subjclass{Primary 20E05}

%

\maketitle

In the paper \cite{bv} the authors give estimates on the growth
rate of the number of \emph{primitive} elements in free groups, as
a function of word length (a primitive element in the free group
of rank $p$ is an element which is a member of a generating set
containing exactly $p$ elements). They produce estimates on the
\emph{exponential growth rate} of the number of such primitive
elements (recall that the exponential growth rate of a set $S$ in
 $F_p$ is defined in  \cite{bv} to be
\begin{equation*}
d_X(S) = \limsup_{N \rightarrow \infty} \left[\frac{\mbox{no. of
elements of $S$ of wordlength at most $N$}} {\mbox{no. of elements
in $F_p$ of wordlength at most $N$}}\right]^{\frac{1}{N}}.
\end{equation*}
The authors produce upper and lower bounds on the exponential
growth rate of the number of primitive elements, but they omit the
following observation:
\begin{theorem}
\label{t2}
 In $F_2$ the exponential growth rate of the number of
primitive elements equals $1/\sqrt{3}.$
\end{theorem}
This observation is significant, since $F_2$ is the only
(nonabelian) free group where the exact value of the exponential
growth rate can be computed.

\begin{proof}[Proof of Theorem \ref{t2}] The result will follow
immediately from  \cite[Proposition 4.1]{bv} and Theorem \ref{t3}
below.
\end{proof}

\begin{theorem}
\label{t3}
 The exponential growth rate of the number of
\emph{cyclically reduced} primitive elements in $F_2$ equals
$1/3.$
\end{theorem}
\begin{proof}
It is sufficient to show that the number of cyclically reduced
primitive elements in $F_2$ grows subexponentially. To show this,
note that each such elements (up to rotation) corresponds to a
\emph{simple geodesic} on the hyperbolic punctured torus (this is
classical, going back to at least Fenchel, but the first reference
known to me is \cite{oz}). Furthermore, the number of simple
geodesics of  (hyperbolic) length bounded above by $L$ on the
punctured torus grows quadratically in $L$ (see \cite{mr2} for
precise results), and for each hyperbolic structure on the
punctured torus there is a constant $C>0,$ such that the ratio of
the hyperbolic length to the word length is at least $1/C$ and at
most $C.$ This is shown in \cite{mr1}. The result follows
immediately.
\end{proof}

\bibliographystyle{amsplain}

\end{document}